\documentclass[12pt]{amsart}
\usepackage{amsmath,amssymb}

\DeclareMathOperator{\Spec}{Spec}

\begin{document}
\newtheorem{thm}{Theorem}[section]
\newtheorem{lem}[thm]{Lemma}
\newtheorem{dfn}[thm]{Definition}
\newtheorem{cor}[thm]{Corollary}
\newtheorem{conj}[thm]{Conjecture}
\newtheorem{clm}[thm]{Claim}
\theoremstyle{remark}
\newtheorem{exm}[thm]{Example}
\newtheorem{rem}[thm]{Remark}
\newtheorem{question}[thm]{Question}
\def\N{{\mathbb N}}
\def\G{{\mathbb G}}
\def\Q{{\mathbb Q}}
\def\R{{\mathbb R}}
\def\C{{\mathbb C}}
\def\P{{\mathbb P}}
\def\Z{{\mathbb Z}}
\def\v{{\mathbf v}}
\def\x{{\mathbf x}}
\def\O{{\mathcal O}}
\def\M{{\mathcal M}}
\def\kbar{{\bar{k}}}
\def\tr{\mbox{Tr}}
\def\id{\mbox{id}}
\def\qed{{\tiny $\clubsuit$ \normalsize}}
\newcommand{\Osh}{\O}
\newcommand{\PP}{\P}
\newcommand{\NN}{\N}
\newcommand{\ZZ}{\Z}
\newcommand{\Hom}{\operatorname{Hom}} 
\newcommand{\ev}{\operatorname{ev}} 
\newcommand{\Hilb}{\operatorname{Hilb}} 
\newcommand{\codim}{\operatorname{codim}} 

\renewcommand{\theenumi}{\alph{enumi}}




\makeatletter
\newcommand*{\da@rightarrow}{\mathchar"0\hexnumber@\symAMSa 4B }
\newcommand*{\da@leftarrow}{\mathchar"0\hexnumber@\symAMSa 4C }
\newcommand*{\xdashrightarrow}[2][]{%
  \mathrel{%
    \mathpalette{\da@xarrow{#1}{#2}{}\da@rightarrow{\,}{}}{}%
   }%
}
\newcommand{\xdashleftarrow}[2][]{%
  \mathrel{%
    \mathpalette{\da@xarrow{#1}{#2}\da@leftarrow{}{}{\,}}{}%
  }%
}
\newcommand*{\da@xarrow}[7]{%
  \sbox0{$\ifx#7\scriptstyle\scriptscriptstyle\else\scriptstyle\fi#5#1#6\m@th$}%
  \sbox2{$\ifx#7\scriptstyle\scriptscriptstyle\else\scriptstyle\fi#5#2#6\m@th$}%
  \sbox4{$#7\dabar@\m@th$}%
  \dimen@=\wd0 %
  \ifdim\wd2 >\dimen@
    \dimen@=\wd2 %
  \fi
  \count@=2 %
  \def\da@bars{\dabar@\dabar@}%
  \@whiledim\count@\wd4<\dimen@\do{%
    \advance\count@\@ne
    \expandafter\def\expandafter\da@bars\expandafter{%
      \da@bars
      \dabar@ 
    }%
  }%
  \mathrel{#3}%
  \mathrel{%
    \mathop{\da@bars}\limits
    \ifx\\#1\\%
    \else
      _{\copy0}%
    \fi
    \ifx\\#2\\%
    \else
      ^{\copy2}%
    \fi
  }%
  \mathrel{#4}%
}
\makeatother

\newcommand{\xmapsto}{\mapstochar\relbar\joinrel\xrightarrow}


\title[Codimension two integral points]{Codimension two integral points on some rationally connected threefolds are potentially dense}

\author{David McKinnon}
\author{Mike Roth}

\begin{abstract}
Let $V$ be a smooth, projective, rationally connected variety, defined over a number field $k$, and let $Z\subset V$ be a closed subset of codimension at least two.  In this paper, for certain choices of $V$, we prove that the set of $Z$-integral points is potentially Zariski dense, in the sense that there is a finite extension $K$ of $k$ such that the set of points $P\in V(K)$ that are $Z$-integral is Zariski dense in $V$.  This gives a positive answer to a question of Hassett and Tschinkel from 2001.
\end{abstract}

\maketitle

\vspace{-5pt}
\section{Introduction}

In \cite{HT}, as Problem 2.13 (``The Arithmetic Puncturing Problem"), Hassett and Tschinkel ask the following question:  

\begin{question}\label{htquestion}
Let $X$ be a projective variety with canonical singularities and $Z$ an algebraic subset of codimension at least $2$, all defined over a number field $k$.  Assume that rational points on $X$ are potentially dense.  Are integral points on $(X,Z)$ potentially dense?
\end{question}

\vspace{.1in}

Of course, the hypothesis that $Z$ has codimension at least two cannot be removed, as there are countless well known examples of varieties with a dense set of rational points but a degenerate set of integral points if $Z$ is a divisor.  Note that when we say that an algebraic set has codimension at least two, we mean that every irreducible component has codimension at least two.

In \cite{HT}, they provide positive answers to this question in various cases, including toric varieties and products of elliptic curves.  The purpose of this paper is to give a positive answer to this question for a large number of examples in dimension up to three.

The case for curves seems vacuous, but if one views a curve defined over a number field as an arithmetic surface, then one can choose $Z$ to be an arithmetic zero-cycle, in which case there is something to prove.  This is Lemma~\ref{curvesintegral}, and is a crucial technical tool for the paper.

For surfaces, the situation is more complicated, as it is unknown which surfaces have a potentially dense set of rational points.  If the Kodaira dimension is negative, however -- which is believed to be the case in which rational points are most plentiful -- we give a positive answer to Question~\ref{htquestion} in Theorem~\ref{kodairanegativesurfaces}.  

In the lengthiest and most difficult part of the paper, $X$ will be a smooth, projective, rationally connected threefold.  It is a theorem of Mori (\cite{Mo}) that there is a birational map $f\colon X\dashrightarrow V$, where $V$ is a normal projective threefold with only $\Q$-factorial and terminal singularities with a morphism $\pi\colon V\to Y$ of one of the following three types:

\begin{enumerate}
\item\label{conic} The variety $Y$ is a normal projective surface with at most rational singularities, and $\pi$ makes $V$ a conic bundle over $Y$.
\item\label{delpezzo} The variety $Y$ is isomorphic to $\P^1$, and a general fibre of the morphism $Y$ is a smooth del Pezzo surface.
\item\label{rankone} The variety $Y$ is a point, and $\mbox{Pic}(V)\cong\Z$.
\end{enumerate}

This list provides a natural set of examples on which to test Question~\ref{htquestion}.  Indeed, in light of Lemma~\ref{blowup}, a positive answer to Question~\ref{htquestion} for the varieties listed above will provide a positive answer for any blowup of such varieties, which constitutes a huge proportion of all smooth, rationally connected threefolds.  In this paper, we will deal with examples from cases \ref{delpezzo} and \ref{rankone}.  Specifically, we prove the following theorems:

\begin{thm}\ref{fano}
Let $X$ be a complex Fano threefold of Picard rank one and index at least two.  Assume that $X$ is defined over a number field $k$, and let $Z$ be an algebraic subset of $X$ of codimension at least two.  If $X$ is a hypersurface of degree 6 in the weighted projective space $\P(1,1,1,2,3)$, then we make the further assumption that $Z$ does not contain the 
unique basepoint $P$ of the square root of the anticanonical linear system.  Then the $Z$-integral points of $X$ are potentially Zariski 
dense.
\end{thm}

\begin{thm}\ref{thm:degreethree}
Let $k$ be a number field.  Let $X$ be a smooth threefold with a map $\pi\colon X\to\P^1$ whose generic fibre is a del Pezzo surface of degree at least three, all defined over $k$.  Let $Z\subset X$ be an algebraic subset of codimension at least two.  Then the $Z$-integral points of $X$ are potentially Zariski dense.
\end{thm}

The rest of the paper is structured as follows.  In section 2, we make some preliminary definitions, and prove two useful lemmas, including Lemma~\ref{curvesintegral}, which is a central technical tool for the rest of the paper.  Section 3 gives a positive answer to Question~\ref{htquestion} for rational and ruled surfaces.  Sections 4 and 6 are the heart of the paper, giving a positive answer to Question~\ref{htquestion} for a wide range of rationally connected threefolds: section 4 deals with Fano threefolds with Picard rank one, and section 6 with del Pezzo fibrations.  In section 5, given a fibration $\pi\colon X\to \P^1$, we prove the existence of sections avoiding a subset of codimension at least two, which may be of independent interest.  Finally, section 7 gives some applications of these results to integral points on families of curves and surfaces, including some classical cases of points integral with respect to a divisor.

\section{Preliminaries}

We first fix some notation and definitions.  Let $X$ be a projective algebraic variety, $Z\subset X$ a Zariski closed subset, both defined over a number field $k$.  Let $M_k$ be the set of places of $k$.  The following definition is essentially Definition 1.4.3 in \cite{Vo}:

\begin{dfn}
Let $S$ be a finite set of places of $k$ containing all the archimedean places.  A subset $R\subset X(k)-Z(k)$ is called {\em $(Z,S)$-integralizable} if and only if there are global Weil functions $\lambda_{Z,v}$ and non-negative real numbers $n_v$ such that $n_v=0$ for all but finitely many places $v$ for each $v$, and such that
\[\lambda_{Z,v}(P)\leq n_v\]
for all $v\in M_k-S$ and $P\in R$.

If the $\lambda_{Z,v}$ and $n_v$ are fixed, we will say that a $k$-rational point $P$ is {\em $(Z,S)$-integral} (or {\em $Z$-integral}, if $S$ is understood) if and only if $\lambda_{Z,v}(P)\leq n_v$ for all $v\in M_k-S$.  
\end{dfn}

This definition is somewhat involved, and for the sake of brevity, we refer the reader to section 1.4 of \cite{Vo} for a more detailed discussion.

We begin with a lemma.  

\begin{lem}\label{blowup}
Let $f\colon X\to Y$ be a birational morphism between irreducible varieties.  Assume that $f$, $X$, and $Y$ are all defined over the same number field $k$.  If Question~\ref{htquestion} has a positive answer for every algebraic subset $Z\subset Y$ of codimension at least two, then it has a positive answer for every algebraic subset $Z\subset X$ of codimension at least two.
\end{lem}

\noindent
{\it Proof:} \/ Let $Z\subset X$ be an algebraic subset of codimension at least two.  Then $f(Z)$ is an algebraic subset of $Y$ of codimension at least two, so by hypothesis the $f(Z)$-integral points of $Y$ are potentially Zariski dense.  But then the $f^{-1}(f(Z))$-integral points of $X$ are potentially Zariski dense as well, so {\em a fortiori} the $Z$-integral points of $X$ are also potentially Zariski dense.  \qed

\vspace{.1in}

The next lemma also appears as Theorem~3.1 in \cite{MZ}, but the proof given there is slightly different.  The idea behind this lemma is to show that if a curve has infinitely many integral points on it, then deleting an arithmetic zero-cycle from it will either delete all the integral points, or else leave an infinite set of integral points.  (If the curve $C$ is projective, then ``integral points'' refers to rational points.)

In the statement of the lemma, $C$ is the curve we're considering, and $Z$ is the ``locus at infinity'' -- that is, we assume that $C$ has an infinite set of $Z$-integral points.  We then delete a further set $N$ which is assumed to intersect $C$ in an arithmetic zero-cycle, and the assumption is that $C$ contains at least one $(Z\cup N)$-integral point.  Lemma~\ref{curvesintegral} then says that $C$ must still contain an infinite set of $(Z\cup N)$-integral points.  
In other words, Question~\ref{htquestion} has a positive answer for curves, considered as arithmetic surfaces.

\begin{lem}\label{curvesintegral}
Let $V$ be an algebraic variety defined over a number field $k$, and let $\mathcal{V}$ be a model of $V$ over $\mbox{Spec}(\O_k)$.  Let $C$ be an irreducible curve on $V$, and let $\mathcal{C}$ be its closure in $\mathcal{V}$.  Let $Z$ and $N$ be Zariski closed subsets of $V$, and let $\mathcal{Z}$ and $\mathcal{N}$ be their closures in $\mathcal{V}$, respectively.  Let $L=Z\cup N$, and let $S$ be a set of places of $k$ that contains all the archimedean places of $k$.

For every place $v$ of $k$ with $v\not\in S$, let $n_v$ be a non-negative real number.  Assume that $n_v=0$ for all but finitely many places $v$.  Choose Weil functions $\lambda_{L,v}$ for each place $v$.  Assume that there is a point $P\in C(k)$ satisfying
\[\lambda_{L,v}(P)\leq n_v\]
for every place $v\not\in S$.

If $N\cap C=\emptyset$ and $C$ contains an infinite set of $(Z,S)$-integral points, then there are infinitely many points $Q\in C(k)$ satisfying
\[\lambda_{L,v}(Q)\leq n_v\]
for every place $v\not\in S$.
\end{lem}

\noindent
{\it Proof:} \/  Since $C(k)$ is infinite, it follows that $C$ must have geometric genus zero or one.  The condition that $C$ contain a dense set of $Z$-integral points implies that $C$ must intersect $Z$ in at most two places of $C$ (places in the sense of points of the normalization of $C$), and that $C\cap Z=\emptyset$ if $C$ has genus one.

We first assume that $C\cap Z=\emptyset$.  Note that without loss of generality, we may assume that $S$ is precisely the set of archimedean places of $k$, as increasing $S$ only makes the lemma easier to prove.  

Let $S'$ be the set of places $v$ of $k$ such that either $v$ is archimedean, or $\mathcal{N}\cap\mathcal{C}$ is supported on $v$.  Note that $S'$ is finite.

For each $v\not\in S'$, $\lambda_{N,v}(Q)=0$ for all $Q\in C$, so we may restrict our attention to $v\in S'$.

If $v$ is finite with corresponding prime $\pi$ of $\O_k$, then the condition $\lambda_{N,v}(P)<n_v$ depends only on the residue class of $P$ modulo a suitable power of $\pi$.  (See for example subsection 2.2.2 of \cite{BG}.)  Thus, the collection of all $Q\in C(k)$ satisfying $\lambda_{N,v}(Q)\leq n_v$ for all finite $v$ contains the set of points $Q$ such that $Q\equiv P\pmod M$ for some suitable nonzero $M\in\O_k$.

There are now two cases: either the geometric genus of $C$ is zero, or one.

If the geometric genus of $C$ is zero, then the theorem follows immediately from the Weak Approximation Theorem for $\P^1$.

Weak Approximation does not hold for curves of genus one, however, so we must work a bit harder.  The set $B$ of points $Q$ such that $Q\equiv P\pmod M$ for some nonzero $M\in \O_k$ contains the image on $C$ of a coset of a finite index subgroup $A$ of the Mordell-Weil group of the normalization $\tilde{C}$ of $C$ over $k$.  Since the set of rational points of $C$ is infinite, the group $A$ is infinite, and so we are done with the case $C\cap N=\emptyset$.

The only cases that remain are when $C$ is a genus zero curve with either one or two places supported on $Z$.  Let $\pi\colon \tilde{C}\to C$ be the normalization map over $k$.  The set of points $R$ of $\tilde{C}$ with $\pi(R)\not\in Z$ are a principal homogeneous space for an arithmetic group ($\mathbb{G}_a$ if there is one place of $C$ on $Z$, and $\mathbb{G}_m$ if there are two places), so by choosing a point $R_0$ on $\tilde{C}$, we can give the $(\pi^*Z,S)$-integral points of $\tilde{C}$ the structure of an arithmetic group $G$.  The set $B$ of points $Q$ of $C$ such that $Q\equiv P\pmod N$ for some nonzero $N\in\O_k$ -- which, as before, is contained in the set of points $Q$ satisfying $\lambda_{L,v}(Q)\leq n_v$ for all $v$ not in $S$ -- contains the image on $C$ of a coset of a finite index subgroup $A$ of $G$, and is therefore infinite, as desired.  \qed

\vspace{.1in}

We will apply Lemma~\ref{curvesintegral} in the case where $Z$ is empty.  

\section{Surfaces}

If Lemma~\ref{curvesintegral} gives a positive answer to Question~\ref{htquestion} for curves, then the next natural question is to ask if it has a positive answer for surfaces.  This is as yet unknown in general, but there are nevertheless a great many cases in which it is known.  For example, Question~\ref{htquestion} has a positive answer for every toric variety, by Corollary 4.2 in \cite{HT}.  In fact, we can prove much more.

\begin{thm}\label{kodairanegativesurfaces}
Let $X$ be a complex surface with negative Kodaira dimension, defined over a number field $k$.  Then Question~\ref{htquestion} has a positive answer for $X$.
\end{thm}

\noindent
{\it Proof:} \/ We begin with the following analogue of Lemma~\ref{curvesintegral}:

\begin{lem}\label{localglobalsurface}
Let $X$ be an algebraic surface, defined over a number field $k$, birational to $\P^2_k$ over $k$.  Fix a model $\mathcal{X}$ for $X$ over $\Spec(\O_k)$, and an effective arithmetic $1$-cycle $\mathcal{Z}$ on $\mathcal{X}$, with $\mathcal{Z}$ defined over $\O_k$.  Let $S$ be a finite set of places of $k$ including all the archimedean places.  For each place $v\not\in S$ of $k$, we fix a non-negative real number $n_v$ and a Weil function $\lambda_{Z,v}$.  If there is a $\mathcal{Z}$-integral point on $X$, then the set of $\mathcal{Z}$-integral points is Zariski dense.
\end{lem}

\noindent
{\it Proof:} \/ Let $f\colon X\dashrightarrow\P^2$ be a birational map defined over $k$.  We may choose $f$ so that there is a finite set of points $P_1\ldots, P_n$ such that $f$ restricts to a birational morphism from $X-\{P_1,\ldots,P_n\}$ to $\P^2$.  In particular, this means that $f(\mathcal{Z})$ has dimension at most $1$ on $\mathcal{X}$.  We may also assume, {\it a fortiori}, that $\mathcal{Z}=f^{-1}(f(\mathcal{Z}))$.  In that case, a point $P\in X$ is $\mathcal{Z}$-integral if and only if $f(P)$ is $f(\mathcal{Z})$-integral.  It therefore suffices to show that the set of $f(\mathcal{Z})$-integral points in $\P^2$ is Zariski dense.  This is well known -- it follows, for example, from Corollary~4.2 of \cite{HT}, or from Theorem~4 of \cite{Sh}.  We include a proof here for completeness.

Let $Q$ be a $\mathcal{Z}$-integral point in $\P^2$, guaranteed by the hypothesis.  Then there is a line $L$ through $Q$ whose closure $\mathcal{L}$ over $\Spec(\Z)$ meets $f(\mathcal{Z})$ in an arithmetic $0$-cycle.  By Lemma~\ref{curvesintegral}, the existence of one $f(\mathcal{Z})$-integral point $Q$ on $L$ implies the existence of infinitely many such points.  For each such point $Q'$, we may find a further line $L'$, different from $L$, that passes through $Q'$ and meets $f(\mathcal{Z})$ in an arithmetic $0$-cycle.  This means that the $f(\mathcal{Z})$-integral points on all the lines $L'$ are also Zariski dense, so the $f(\mathcal{Z})$-integral points of $\P^2$ -- and therefore also of $X$ -- are Zariski dense, as desired.  \qed

\vspace{.1in}

Every rational surface is well known to be the blowup of some Hirzebruch surface or of the projective plane.  Since all of those are toric varieties, Question~\ref{htquestion} has a positive answer for them.  Therefore, by Lemma~\ref{blowup}, Question~\ref{htquestion} has a positive answer for every rational surface.

If $X$ is a ruled surface, then there is a fibration $f\colon X\to C$ for some smooth curve $C$.  If $C$ has genus at least two, then Question~\ref{htquestion} has a vacuously positive answer for $X$, because the rational points on $X$ are not potentially dense.  If the genus of $C$ is zero, then $X$ is rational and Question~\ref{htquestion} has a non-vacuously positive answer for $X$, as just noted.  

Thus, assume that $C$ has genus 1, and let $Z$ be an algebraic subset of codimension at least two -- that is, let $Z$ be a finite set of points of $X$.  By a finite extension of the field of definition $k$, we may assume that there is a $Z$-integral point $P$ on $X$ defined over $k$, and that $C$ has an infinite number of $k$-rational points.  By, for example, Theorem V.2.17.(c) of \cite{Ha}, there is a very ample divisor class $V$ on $X$ whose elements are sections of $f$, and therefore have infinitely many rational points.  Let $Y_1$ be a curve in the class $V$ that contains the point $P$, but does not intersect $Z$.  By Lemma~\ref{curvesintegral}, $Y_1$ has a Zariski dense set of $Z$-integral points, and in particular, there are an infinite number of fibres $F$ of $f$ for which $Y_1\cap F$ is a $Z$-integral point, and for which $F\cap Z=\emptyset$.  By Lemma~\ref{curvesintegral} again, this means that $F$ has a dense set of $Z$-integral points, implying that the set of $Z$-integral points is dense, and that Question~\ref{htquestion} has a positive answer for $X$.

Every surface of negative Kodaira dimension is the blowup of a rational or ruled surface.  Thus, by Lemma~\ref{blowup}, Theorem~\ref{kodairanegativesurfaces} is proven.  \qed

\vspace{.1in}

For surfaces with non-negative Kodaira dimension, the situation is more complex, and indeed it is still not known which of these surfaces have a Zariski dense set of rational points, never mind integral ones.  We will therefore move on to threefolds.

\section{Fano threefolds}

For the purposes of this paper, a Fano threefold is a smooth, three-dimensional algebraic variety $X$ whose anticanonical sheaf $-K_X$ is ample.  If $X$ has Picard rank one -- that is, if the Picard group of $X$ is isomorphic to $\Z$ -- then there is a unique ample generator $H$ of the Picard group.  The index of a Fano threefold is the unique integer $r$ such that $-K_X=rH$.  The main theorem of this section is the following:

\begin{thm}\label{fano}
Let $X$ be a complex Fano threefold of Picard rank one and index at least two.  Assume that $X$ is defined over a number field $k$, and let $Z$ be an algebraic subset of $X$ of codimension at least two.  If $X$ is a hypersurface of degree 6 in the weighted projective space $\P(1,1,2,3)$, then we make the further assumption that $Z$ does not contain the unique basepoint $P$ of the square root of the anticanonical linear system.  Then the $Z$-integral points of $X$ are potentially Zariski dense.
\end{thm}

\noindent
{\it Proof:} \/ The proof relies crucially on the classification of Fano threefolds of Picard rank one, found (for example) in \cite{IP}.  Section 12.2 of \cite{IP} gives the following list of Fano threefolds of Picard rank one and index at least two:

\begin{enumerate}
\item $\P^3$
\item A smooth quadric in $\P^4$.
\item A smooth linear section of the Pl\"ucker-embedded Grassmannian $\mbox{Gr}(2,5)$.
\item A smooth intersection of two quadrics in $\P^5$.
\item A smooth cubic hypersurface in $\P^4$.
\item A double cover of $\P^3$, branched on a smooth quartic surface.
\item A smooth hypersurface of degree 6 in the weighted projective space $\P(1,1,1,2,3)$.
\end{enumerate}

Note that our list is in the reverse order of that in \cite{IP}, and in particular, item (g) is indeed the same as the unnamed threefold with $-K_X^3=8$ and $h^{1,2}=21$.  

We now proceed by cases.  Case (a) is the easiest, as the answer to question~\ref{htquestion} is well known to be positive for $\P^3$ -- see for example \cite{HT} or \cite{Sh}.

\vspace{.1in}

\noindent
{\it Case (b):} \/ $X$ is a smooth quadric hypersurface in $\P^4$.

\vspace{.1in}

Let $P$ be a $k$-rational point of $X-Z$, and let $\pi\colon X\to\P^3$ be the linear projection from $P$.  Then $\pi(Z)$ is an algebraic subset of $\P^3$ of codimension at least two, and so the $\pi(Z)$-integral points are potentially Zariski dense.  It therefore follows immediately that the $Z$-integral points of $X$ are also potentially Zariski dense.

\vspace{.1in}

\noindent
{\it Case (c):} \/ $X$ is a smooth linear section of the Pl\"ucker-embedded Grassmannian $\mbox{Gr}(2,5)$.

\vspace{.1in}

$X$ can be obtained by blowing up a smooth quadric threefold $Q\subset\P^4$ along a smooth rational curve of degree three, and then contracting the strict transform of a smooth quadric surface.  If $\pi\colon Y\to Q$ is the blowup, and $\phi\colon Y\to X$ is the contraction, then a $Z$-integral point on $X$ corresponds to a $\phi^*Z$-integral point on $Y$.  Any $\pi(\phi^*Z)$-integral point of $Q$ pulls back to a $\phi^*Z$-integral point of $Y$, so it suffices to show that the $\pi(\phi^*Z)$-integral points of $Q$ are potentially Zariski dense.

The scheme $\pi(\phi^*Z)$ is contained in the union of a smooth quadric surface and a subset of $Q$ of codimension at least two.  It therefore suffices to prove that the $Z$-integral points of $Q$ are potentially Zariski dense, where $Z$ is the union of a smooth quadric surface $S$ and a subset $W$ of codimension at least two.

After a finite extension of the base field $k$, we may assume that there is a $Z$-integral point $P$ on $Q$, and that the group $\O_k^*$ is infinite, where $\O_k^*$ is the group of units of the ring of integers $\O_k$ of $k$.  Let $T$ be a $2$-plane containing $P$, but with $T\cap S$ finite, $T\cap Q$ irreducible and smooth at $P$, and $T\cap W=\emptyset$.  Then $(T\cap Q)-(T\cap S)$ is a rational curve with at most two places deleted, and has a $Z$-integral point.  Therefore, by Lemma~\ref{curvesintegral}, $T\cap Q$ contains infinitely many $Z$-integral points.

For each such point $P'$, we can find another $2$-plane $T'$ such that $P'\in T'$, $T'\cap S$ finite, $T'\cap Q$ irreducible and smooth at $P'$, $T'\cap W=\emptyset$, and $T'\neq T$.  This, via Lemma~\ref{curvesintegral}, yields a set of $Z$-integral points whose Zariski closure $Y$ has dimension at least $2$.  If $Y\neq Q$, then for each $Z$-integral point $P_i''$ on $Y$, we can find a $2$-plane $T_i''$ such that $P_i''\in T_i''$, $T_i''\cap S$ finite, $T_i''\cap Q$ irreducible and smooth at $P_i''$, $T_i''\cap W=\emptyset$, and $T_i''\cap Q\not\in Y\cup T_1''\cup\ldots\cup T_{i-1}''$.  By Lemma~\ref{curvesintegral}, we obtain a set of $Z$-integral points that is Zariski dense in $Q$.

\vspace{.1in}

\noindent
{\it Case (d):} \/ $X$ is a smooth intersection of two quadrics in $\P^5$.

\vspace{.1in}

After a finite extension of the base field $k$, we can choose a $Z$-integral point $P$ that is not contained in $Z$, and such that the singular locus of the linear projection of $X$ from $P$ is not contained in the image of $Z$.  Let $\pi_1\colon X\to\P^4$ be the projection from $P$.  Then $\pi_1(X)$ is a singular cubic threefold, and if $P'$ is a singular point of $\pi_1(X)$ that is not contained in $\pi_1(Z)$, then the projection $\pi_2\colon\pi_1(X)\to\P^3$ of $\pi_1(X)$ away from $P'$ induces a birational map $\phi\colon X\to\P^3$ such that $\phi(Z)$ is an algebraic subset of $\P^3$ of codimension at least two.  Since $\phi(X)$-integral points are well known to be potentially Zariski dense in $\P^3$, it follows that $Z$-integral points on $X$ are also potentially Zariski dense, as desired.

\vspace{.1in}

\noindent
{\it Case (e):} \/ $X$ is a smooth cubic threefold in $\P^4$. 

\vspace{.1in}

Let $\ell$ be a line on $X$ with $\ell\cap Z=\emptyset$, and let $\pi\colon Y\to X$ be the blowing up of $X$ along $\ell$.  Then $Y$ admits the structure of a conic bundle $\phi\colon Y\to\P^2$, where the exceptional divisor $S$ of $\pi$ is a rational surface and a double section of $\phi$.

After a fixed extension of $k$, we may assume that $\ell$, $Y$, $\pi$, and $\phi$ are all defined over $k$, and that $S$ has a dense set of $k$-rational points, including a point $P$ that is also $Z$-integral.  By Lemma~\ref{localglobalsurface}, this means that $S$ has a Zariski dense set of $Z$-integral points as well.

The dimension of $Z$ is at most one, so there is a dense set $A$ of points of $S$ such that for all $Q\in A$, the fibre of $\phi$ through $Q$ does not meet $Z$.  by Lemma~\ref{curvesintegral}, each such fibre has an infinite set of $Z$-integral points, and so the $Z$-integral points of $X$ are dense.

\vspace{.1in}

\noindent
{\it Case (f):} \/ $X$ is a double cover of $\P^3$ branched on a smooth quartic surface.

\vspace{.1in}

The threefold $X$ is known to be unirational (see for example \cite{IP}, Example 10.1.3.(iii)), so we may extend the number field $k$ to ensure that $X$ has a Zariski dense set $S$ of rational points.  Let $\pi\colon X\to \P^3$ be the double cover.  The set $\pi(S)$ is Zariski dense in $\P^3$, and by extending the field $k$ again we may assume that at least one point $P$ of $\pi(S)$ is $\pi(Z)$-integral.

Any line $\ell$ in $\P^3$ lifts to a curve of geometric genus at most one on $X$.  The net $N$ of lines through $P$ induces an elliptic threefold structure (fibred over a rational surface) on a blowup $\tilde{X}$ of $X$.  By Merel's Theorem on the uniform boundedness of torsion on elliptic curves (see \cite{Me}), there is a proper Zariski closed subset $G$ of $\tilde{X}$  which contains all the $k$-rational points of $\tilde{X}$ that are torsion points on their fibre.  We further enlarge $G$ to contain all the singular fibres of $\tilde{X}$.  Let $T$ be the complement of the image of $G$ in $\P^3$, intersected with the set $\pi(S)$.  Then $T$ consists entirely of $k$-rational points $Q$ of $\P^3$ whose preimages on $X$ are also $k$-rational, and such that the elliptic curve lying over the line joining the point $Q$ to $P$ has positive rank.  (The point $P$ is viewed as the identity element.)

Since $Z$ has codimension at least two, there is a Zariski dense set of points $Q$ in $T$ such that the line $\ell$ joining $P$ to $Q$ is disjoint from $Z$.  In each such case, the elliptic curve $E$ lying over $\ell$ has positive Mordell-Weil rank (because $Q$ is non-torsion with respect to $P$), and so by Lemma~\ref{curvesintegral}, $E$ contains an infinite set of $Z$-integral points.  Since the set of such $E$ is Zariski dense, the theorem follows.

\vspace{.1in}

\noindent
{\it Case (g):} \/ $X$ is a smooth hypersurface of degree $6$ in the weighted projective space $\P(1,1,1,2,3)$, and the subset $Z$ does not contain the unique basepoint $P$ of the square root of the anticanonical linear system.

\vspace{.1in}

In this case, $X$ is a double cover of the cone $V$ in $\P^6$ over the Veronese embedding of $\P^2$ in $\P^5$; denote the cover by $\pi\colon X\to V$.  Note that $P$ is the preimage $P=\pi^{-1}(v)$ of the vertex $v$ of the cone $V$.  

Blowing up $v$ on $V$ yields a smooth threefold $V^*$, which admits the structure of a $\P^1$-bundle $f\colon V^*\to\P^2$ over $\P^2$, where the fibres of the bundle are the strict transforms of the lines of the ruling of $V$.  The corresponding blowup of $X$ yields a biregular map $g\colon X^*\to X$ and a double cover $\pi^*\colon X^*\to V^*$, where $X^*$ inherits the structure of an elliptic fibration over $\P^2$ via $\phi=f\circ\pi^*$.

Since $P\not\in Z$, it follows that $Z^*=g^{-1}(Z)$ is at least codimension two as a subset of $X^*$.  In \cite{BT}, the authors show that there is a two-dimensional family of double sections of $\phi$ that are singular, but birational to $K3$ surfaces.  After a possible finite field extension, we may assume that one of those double sections, which we will call $S$, satisfies the following properties:
\begin{itemize}
\item $S$ intersects $Z$ properly.
\item $S$ contains a singular point $s$ which is $Z^*$-integral.
\item $S$ contains a Zariski dense set of rational points.
\end{itemize}

To see that such a choice is possible, note that \cite{BT} proves the Zariski density of the rational points, and allows for a two-dimensional linear system full of such double sections $S$.  (The extra field extension is necessary for the existence of the $Z^*$-integral singular point.)

Given such an $S$, we blow up the singular locus with $h\colon S^*\to S$ to obtain an elliptically fibred, smooth $K3$ surface $S^*$.  The exceptional divisor over $s$ is a $(-2)$-curve on $S^*$ with a dense set of rational points, each of which is $h^{-1}(Z^*)$-integral.  Thus, every smooth elliptic curve in the elliptic fibration on $S^*$ contains at least one $h^{-1}(Z^*)$-integral point.  The density of rational points on $S^*$ implies that there are infinitely many such curves with positive Mordell-Weil rank.  Therefore, since $Z^*\cap S$ is of codimension at least two, we conclude by Lemma~\ref{curvesintegral} that the set of $h^{-1}(Z^*)$-integral points on $S^*$ is Zariski dense, and therefore that the $Z^*$-integral points on $S$ are also Zariski dense on $S$.

For any $Z^*$-integral point $x$ on $X^*$, Lemma~\ref{curvesintegral} again shows that the $Z^*$-integral points are Zariski dense on the fibre of $\phi$ through $x$, provided that its Mordell-Weil rank is positive.  Since \cite{BT} proves that the set of rational points on $X^*$ are Zariski dense, there is a Zariski dense set of fibres with positive Mordell-Weil rank.  Therefore, the set of $Z^*$-integral points on $X^*$ is Zariski dense.  This immediately implies that the set of $Z$-integral points on $X$ is Zariski dense, as desired.  \qed

\section{Sections avoiding given subsets}

In this section we prove a lemma guaranteeing the existence of a section of a rationally connected fibration over a curve, such that the section avoids (respectively fails to be contained in) a given subset of codimension $\geq 2$  (respectively $\geq 1$).  The result is well-known to experts on families of curves on varieties, but we include a proof for lack of a reference.   We begin by recalling background material.

Let $X$ be a smooth projective variety, and set $n=\dim(X)$. 

Recall that a {\em rational curve} in $X$ is a nonconstant map $f\colon\PP^1\longrightarrow X$.  The rational curve
is said to be free if $f^{*}T_{X} = \oplus_{i=1}^{n} \Osh_{\PP^1}(a_i)$, with each $a_i\geq 0$.
Fix an ample line bundle $L$ on $X$. Then for each $d\geq 1$ 
there is a quasi-projective variety $\Hom(\PP^1,X)_d$ parameterizing
maps $f\colon \PP^1\longrightarrow X$ such that $\deg(f^{*}L)=d$. (This is a special case of the general
construction of \cite[1.10]{Ko} constructing parameter spaces $\Hom(Y,X)$ for any projective varieties $Y$ and $X$.  
Any morphism $Y\longrightarrow X$ can be identified with its graph, a subset of $Y\times X$, and the spaces
$\Hom(Y,X)$ are then realized as the open subscheme of $\Hilb(Y\times X)$ parameterizing such graphs.
The restriction $\deg(f^{*}L)=d$ is used to fix the Hilbert polynomial of the graph. )

For a map $f\colon \PP^1\longrightarrow X$, with $\deg(f^{*}L)=d$, we denote by $[f]$ the corresponding 
point of $\Hom(\PP^1,X)_d$.  
One also has an {\em evaluation map} 
$$\begin{array}{rcccc}
\ev & \colon & \Hom(\PP^1,X)_{d}\times\PP^1 & \xrightarrow{\rule{0.75cm}{0cm}} & X \\
& & ([f],p) & \xmapsto{\rule{1cm}{0cm}} & f(p). \rule{0cm}{0.6cm}\\
\end{array}
$$

Let $\Hom(\PP^1,X)^{\circ}_{d}$ denote the subset of $\Hom(\PP^1,X)_{d}$ 
consisting of those $[f]$ such that $f$ is free.  By \cite[II.3.5.4, p. 115]{Ko}, $\Hom(\PP^1,X)_{d}^{\circ}$ 
is an open subset of $\Hom(\PP^1,X)_{d}$, and the evaluation map 
$$\Hom(\PP^1,X)_{d}^{\circ}\times \PP^1 \xrightarrow{\rule{0.3cm}{0cm}\ev\rule{0.3cm}{0cm}} X$$
is smooth. (Thus $\Hom(\PP^1,X)_{d}^{\circ}$ is also smooth, although one can see this last point directly by 
computing the tangent space to the Hilbert scheme).

\begin{lem}\label{lem:avoiding-subsets} (a) Let $X$ be a smooth irreducible projective variety defined over an 
algebraically closed field of characteristic zero, $\pi\colon X\longrightarrow \P^1$ a
surjective map whose general fibre is rationally connected, $Z\subset X$ a subvariety of codimension $\geq 2$, 
and $T\subset X$ a subvariety
of codimension $\geq 1$.  Then there exists a section of $\pi$ which is not contained in $T$, and which does not
meet $Z$. (b) If $X$, and $\pi$, $Z$, and $T$ are defined over a field $k$ of characteristic zero, and if the
general fibre of $\pi$ over $\overline{k}$ 
is rationally connected, then there exists such a section defined over a finite extension $k'$ of $k$. 
\end{lem}

\noindent
{\it Proof:} \/ 
We first prove (a). 
By \cite[Theorem 1.1]{GHS} there is a map $g\colon \PP^1\longrightarrow X$ which is a section of $\pi$.
Furthermore, by \cite[2.13]{KMM}, given that such a section exists, and given any point $q$ on a smooth fibre
of $\pi$, there exists a curve $f'\colon \PP^1 \longrightarrow X$ which is a free curve, a section of $\pi$, and
passes through $q$ (i..e, so that $q$ is in the image of $f'$).

Choose any point $q$ in a smooth fibre, and not in $Z$ or $T$, and let $f'$ be a free curve and section 
passing through $q$ provided by those theorems.  Set $d=\deg((f')^{*}L)$, and let $V$ be the irreducible component
of $\Hom(\PP^1,X)_{d}^{\circ}$ containing $[f']$.  
We consider the diagram

$$
\begin{array}{ccc}
V\times \PP^1 & \xrightarrow{\rule{0.3cm}{0cm}\ev\rule{0.3cm}{0cm}} & X \\
\phantom{\scriptsize p_1}\Bigg\downarrow {\scriptsize p_1} \\
V \\
\end{array}.
$$
The property that a rational curve $f$ is a section of $\pi$ is equivalent to $\deg(f^{*}\pi^{*}\Osh_{\PP^1}(1))=1$,
i.e., that the degree of $\ev^{*}\pi^{*}\Osh_{\PP^1}(1)$ on the fibre $p_1^{-1}([f])$ is $1$.  Since
$p_1$ is flat, the degree of $\ev^{*}\pi^{*}\Osh_{\PP^1}(1)$ is constant on the fibres of $p_1$ and it follows
that every $[f]\in V$ is also a section of $\pi$.

The map $p_1$ is proper, and by \cite[II.3.5.4.2, p. 115]{Ko}, $\ev$ is smooth.  
The set of $[f]\in V$ such that $f(\PP^1)$ is 
contained in $T$ is the locus where the map $\ev^{-1}(T)\stackrel{p_1}{\longrightarrow} V$ has $1$-dimensional
fibres.  By upper semicontinuity of fibre dimension this locus is a closed subset of $V$.  Let $U'$ be its
complement.  The set $U'$ is nonempty since $[f']\in U'$.   Every point $[f]\in U'$ is now a section of $\pi$
not contained in $T$.  To prove part (a) we just need to find such an $[f]$ so that $f(\PP^1)\cap Z=\emptyset$. 

Set $N=\dim(V)$. Since $Z$ is of codimension $\geq 2$, and $\ev$ smooth, $\ev^{-1}(Z)$ also has codimension $\geq 2$,
and hence has dimension at most $N+1-2=N-1$.  Thus $p_1(\ev^{-1}(Z))$ has dimension $\leq N-1$ and so is a proper 
subset of $V$.  Let $U''$ be its complement.  Any $[f]\in U''$ satisfies $f(\PP^1)\cap Z=\emptyset$.  Since $V$
is irreducible, $U:=U'\cap U''\neq \emptyset$, proving (a).

To see (b) we first note that if $X$ is defined over $k$, then we can choose an ample $L$ defined over $k$, and 
then $\Hom(\PP^1,X)_{d}$ is also defined over $k$ for each $d\in \NN$. (As above, one starts with the Hilbert
scheme $\Hilb(\PP^1\times X)$ and restricts to the open subset which are the graphs of morphisms.  The Hilbert
scheme and the condition of being a graph can be expressed over $k$.)  The open condition that a morphism $f$
is free is similarly defined over $k$, as is the condition that $f$ is a section of $\pi$ (this being again
a condition on the degree of $\ev^{*}\pi^{*}\Osh_{\PP^1}(1)$ on the fibres of $p_1$).

If $T$ and $Z$ are defined over $k$, then the closed locus where the fibre dimension of 
$\ev^{-1}(T)\longrightarrow \Hom(\PP^1,X)_{d}^{\circ}$ is $1$ is defined over $k$, and so is the 
closed subset $p_1(\ev^{-1}(Z))$.  Thus, the intersection of their complements is also defined over $k$.
For each $d\in \NN$ 
we let $U_d\subseteq\Hom(\PP^1,X)_{d}$ be the intersection of the complements, 
along with the intersection with the open conditions of the maps being free, and being a section of $\pi$.

By part (a), over $\overline{k}$, there is some
$d$ for which the corresponding $U_d$ is nonempty.  But, if $U_d$ is nonempty after base extension, 
it was nonempty to begin with.  
Since $\Hom(\PP^1,X)_{d}$ is of finite type, the residue field of any closed point is finite over $k$. Thus
taking any closed point $[f]\in U_d$, and letting $k'$ be its residue field, we obtain a section 
$f\colon \PP^1\longrightarrow X$ defined over $k'$, avoiding $Z$, and not contained in $T$. \qed

Using an idea from \cite{GHS} due to de Jong, one can extend 
Lemma \ref{lem:avoiding-subsets} to the case where the base curve has
arbitrary genus.   We will not need this extension, but record the statement and the idea of its proof.

\begin{cor}
(a) Let $X$ be a smooth irreducible projective variety defined over an algebraically
closed field of characteristic zero, $\pi\colon X\longrightarrow C$ a
surjective map to a smooth curve $C$ such that the general fibre of $\pi$ is rationally connected, 
$Z\subset X$ a subvariety of codimension $\geq 2$, and $T\subset X$ a subvariety
of codimension $\geq 1$.  Then there exists a section of $\pi$ which is not contained in $T$, and which does not
meet $Z$. (b) If $X$, and $\pi$, $Z$, and $T$ are defined over a field $k$ of characteristic zero, and if the
general fibre of $\pi$ over $\overline{k}$ 
is rationally connected, then there exists such a section defined over a finite extension $k'$ of $k$. 
\end{cor}

{\it Proof:} \/ 
We repeat the argument of de Jong from \cite[\S3.2]{GHS}.  To prove (a), given $\pi\colon X\longrightarrow C$ 
choose any finite map $g\colon C\longrightarrow \PP^1$, and then form the ``norm'' of $X$.  This is a variety and
map $\varphi\colon Y\longrightarrow\PP^1$ (well defined up to birational equivalence) whose fibre over a general 
point $p\in \PP^1$ is the product $\prod_{q\in g^{-1}(p)} \pi^{-1}(q)$.  
The utility of the norm construction is that sections of $\varphi$ give sections of $\pi$.
Given a section $\sigma$ of $\varphi$, for each $p\in \PP^1$, $\sigma$ gives a point of 
$\prod_{q\in g^{-1}(p)} \pi^{-1}(q)$, and thus for each point $q\in C$, setting $p=g(q)$, $\sigma$ gives a
point in the fibre $\pi^{-1}(q)$. 

To ensure that the resulting section of $\pi$ misses $Z$ and is not contained in $T$, we define appropriate
subsets of $Y$.  
Let $\tilde{Z}\subset Y$ be the subset 
$$\tilde{Z} = \left\{ y\in Y \,\,\left|\,\, 
\mbox{\begin{minipage}{0.55\textwidth}
at least one of the coordinates of  $y \in \varphi^{-1}(\varphi(y)) =
\prod_{q\in g^{-1}(\varphi(y))} \pi^{-1}(q)$ is in $Z$
\end{minipage}
}\right.\right\}$$
and similarly define $\tilde{T}$.

Sections $\sigma$ of $\varphi$ which do not meet $\tilde{Z}$ and are not contained in $\tilde{T}$ induce 
sections of $\pi$ similarly missing $Z$ and not contained in $T$.  The codimensions of $\tilde{Z}$ in $Y$
is equal to the codimension of $Z$ in $X$, and similarly $\codim(\tilde{T},Y)=\codim(T,X)$. 

Since the product of rationally connected
varieties is rationally connected, the general fibre of $Y$ is rationally connected, and so we can apply 
Lemma \ref{lem:avoiding-subsets}(a), proving (a) of the corollary.

To prove (b), supposing everything defined over $k$, if we choose our map $g\colon C\longrightarrow \PP^1$
to be defined over $k$, then so are $Y$, $\tilde{Z}$, and $\tilde{T}$.  Thus applying
Lemma \ref{lem:avoiding-subsets}(b),  we obtain a section of $\varphi$ defined over a finite extension $k'$
missing $\tilde{Z}$ and not contained in $\tilde{T}$. This then induces a section of $\pi$, also defined over $k'$,
with the desired properties.  \qed

\section{Del Pezzo fibrations}

in this section, we prove the potential density of integral points for del Pezzo fibrations, provided that the degree of the (generic) del Pezzo surface is at least three.

Let $\pi\colon X\to Y$ be a morphism, where $X$ is a smooth, rationally connected, projective threefold, and $Y$ is a smooth curve.  Since $X$ is rationally connected, $Y$ must be isomorphic to $\P^1$ over $k$.  (This may require a finite extension of $k$.)  We further assume that a general fibre of $\pi$ is a del Pezzo surface.  

Choose models $\mathcal{X}$ and $\mathcal{Y}$ for $X$ and $Y$, respectively, over $\O_k$, and extend $\pi$ to a rational map from $\mathcal{X}$ to $\mathcal{Y}$.  Let $\mathcal{Z}\subset\mathcal{X}$ be a closed subscheme of codimension at least two.  We will show that in many cases, the $\mathcal{Z}$-integral points of $\mathcal{X}$ are potentially Zariski dense.

\begin{thm}\label{thm:degreethree}
Let $k$ be a number field.  Let $X$ be a smooth threefold with a map $\pi\colon X\to\P^1$ whose generic fibre is a del Pezzo surface of degree at least three, all defined over $k$.  Let $Z\subset X$ be an algebraic subset of codimension at least two.  Then the $Z$-integral points of $X$ are potentially Zariski dense.
\end{thm}

\noindent
{\it Proof:} \/ 
Let $T$ be the union of the $(-1)$-curves in the fibres of $\pi$.  Applying Lemma \ref{lem:avoiding-subsets}, after at most a finite field extension -- which we continue to call $k$ -- we obtain a $k$-rational section $\sigma\colon\P^1\to X$ of $\pi$ whose image is a smooth rational curve $C\subset X$,  and disjoint from $Z$, and meeting $T$ in only finitely many points (i.e., only finitely many points of $C$ are contained in $(-1)$-curves of the fibres of $\pi$).  Furthermore, after blowing up, we may decrease the degree of the generic fibre of $\pi$ to three without changing the hypothesis or conclusion of the Theorem~\ref{thm:degreethree}.  (We choose the blowup locus to be disjoint from $Z$.)  Let $S\subset \PP^1(k)$ be the finite subset of points $p$ where either $\pi^{-1}(p)$ contains a $1$-dimensional  component of $Z$, or $\pi^{-1}(p)$ intersects $C$ in a point on a $(-1)$-curve of the fibre.

Theorem~\ref{thm:degreethree} then follows by applying the following lemma to the fibres $\pi^{-1}(p)$, with $p\in \PP^1(k)\setminus S$.  (Note that Lemma~\ref{delpezzodense} is not implied by Lemma~\ref{localglobalsurface} because a del Pezzo surface need not be birational to $\P^2$ over $k$.)

\begin{lem}\label{delpezzodense}
Let $V$ be a del Pezzo surface of degree three defined over a number field $k$, and let $\mathcal{V}$ be a model for $V$ over $\Spec(\O_k)$.  Let $\mathcal{Z}\subset\mathcal{V}$ be an algebraic subset of codimension at least two.  Assume that there is a $k$-rational point $P\in V(k)$ that is $\mathcal{Z}$-integral, and that does not lie on a $(-1)$-curve of $V$.  Then the $\mathcal{Z}$-integral points are Zariski dense.
\end{lem}

\noindent
{\it Proof of lemma:} \/ 
A general member of the linear system $|-K_V|$ is a smooth curve of genus one, and $|-K_V|$ is basepoint free because $V$ is del Pezzo.  Consider the linear subsystem of $|-K_V|$ consisting of curves containing $P$.  It has dimension three, so we can choose a pencil $H$ of curves defined over $k$ such that every curve in $H$ contains $P$, a general curve in $H$ is smooth, and the base locus consists of three points $\{P,Q,R\}$, none of which lie in $Z$.  (Note that $P\not\in Z$ trivially.)

Since $H$ is defined over $k$, so is the triple $\{P,Q,R\}$, and we may therefore blow it up to obtain a surface $\tilde{V}$ defined over $k$, with a morphism $\psi\colon\tilde{V}\to V$ whose fibres are precisely the (strict transforms of) the curves in $H$.  The morphism $\psi$ makes $\tilde{V}$ into an elliptic surface, with a section $\O$ given by the exceptional curve lying over $P$.  Note that $\O$ is disjoint from $Z$.  Since $\psi$ has a section, it has no multiple fibres, so we may enlarge $Z$ to contain all the singular points of fibres of $\psi$.

The class $-K_V$ embeds $V$ in $\P^3$ as a smooth cubic surface.  We may therefore consider the curve $T$ defined by the intersection of the tangent plane $T_P$ with the embedded surface $V$.  Note that $T$ is irreducible because $P$ does not lie on any $(-1)$-curves, so $T$ is an irreducible plane cubic curve.  Moreover, $T$ is singular at $P$, so it has geometric genus zero.  Indeed, $T$ is birational to $\P^1$ over $k$ via projection from $P$ in the plane $T_P$, so $T$ has a dense set of rational points.

For each rational point $A$ of $T$, the intersection of $T_A$ with the embedded surface $V$ is again a cubic curve with a singularity at $A$, albeit possibly reducible.  At most finitely many $A$ correspond to reducible curves in this way (there are only finitely many intersections of $T$ with lines), so there are infinitely many $A$ whose tangent curves are birational to $\P^1$ over $k$, and therefore have a dense set of rational points.  We therefore deduce that the rational points of $V$ are Zariski dense.

This means that there is a dense set of $k$-rational points on $V$ each lying on a smooth fibre of $\psi$ and having infinite order in that fibre.  (To see this, note that by a theorem of Merel (\cite{Me}), there is a positive integer $N$ such that for any elliptic curve defined over $k$, and any $k$-rational point $A$ of finite order, the order of $A$ divides $N$.  Therefore, the set of $k$-rational points of finite order in their fibre is not Zariski dense in $V$.)  In particular, there are an infinite number of genus one curves on $V$ that contain an infinite set of $k$-rational points, one of which is the $Z$-integral point $P$.  By Lemma~\ref{curvesintegral}, each of those curves contains an infinite set of $Z$-integral points.  We conclude that the $Z$-integral points on $V$ are dense.  \qed

\vspace{.1in}

We now finish the proof of Theorem~\ref{thm:degreethree}.  The curve $C$ is disjoint from $Z$ over the generic fibre, so after a further finite extension of $k$ -- which we stubbornly persist in calling $k$ -- we may assume that $C$ contains a $Z$-integral point.  To see this, note that over $\mbox{Spec}(\O_k)$, $\mathcal{N}=\mathcal{C}\cap\mathcal{Z}$ is an arithmetic zero-cycle supported on finitely many places of $k$.  After a suitably chosen finite extension of $k$, we may assume that for every place $v$ of $k$, there is a point $p_v$ of $\mathcal{C}$ lying over $v$ that does not lie in the support of $\mathcal{N}$.  By the Chinese Remainder theorem -- since $C$ is a rational curve -- there is some $k$-rational point $P$ of $C$ such that for all $v$ over which $\mathcal{N}$ is supported, $P\equiv p_v\mod v$.  This $P$ is the $Z$-integral point that we seek.

Since $C$ is a rational curve, Lemma~\ref{curvesintegral} implies that $C$ contains an infinite number of $Z$-integral points.  For each such point, the corresponding fibre, by Lemma~\ref{localglobalsurface}, contains a dense set of $Z$-integral points.  It therefore follows that the $Z$-integral points of $V$ are Zariski dense, as desired.  \qed

\section{Application to integral points in families}

We can apply the theorems of the previous sections to families of curves and surfaces on surfaces and threefolds, to get results about integral points in the classical sense -- that is, integral points with respect to a divisor.  The proof of the following theorem is trivial:

\begin{thm}\label{dPsurfaceintegral}
	Let $X$ be a smooth, projective variety of dimension $n$, defined over a number field $k$, and let $\mathcal{P}$ be an $(n-m)$-dimensional linear system of $m$-dimensional cycles on $X$ whose base locus $Z$ is of pure dimension $m-1$.  (In other words, elements of a basis for $\mathcal{P}$ intersect properly.)  If the $Z$-integral points of $X$ are Zariski dense, then there is a Zariski dense set of cycles in $\mathcal{P}$ with at least one $Z$-integral point defined over $k$. 
\end{thm}

Note that $Z$ is a divisor on each of the cycles in $\mathcal{P}$, so that a $Z$-integral point on a surface $S$ in the family is an integral point in the classical sense.  Indeed, if the divisors in $\mathcal{P}$ are ample on $X$, then the base locus $Z$ is the intersection of any $m+1$ distinct cycles in $\mathcal{P}$, so $Z$ is an ample divisor on any smooth, irreducible cycle in the system.

The results of this paper show that the hypotheses of Theorem~\ref{dPsurfaceintegral} are satisfied when $X$ is any rational or ruled surface, or any rationally connected threefold of a type considered in sections 4 or 5.

In particular, though, the case $m=1$, in which the cycles in $\mathcal{P}$ are curves, is particularly interesting.  For any subset $Z'\subset Z$, any $Z$-integral point is automatically also $Z'$-integral.  Therefore, when $m=1$, the reduced induced subscheme of $Z$ is a finite set of points, so we may assume that $Z$ is a single point, or even a different single point of $Z$ for every curve in $\mathcal{P}$, and the conclusion of Theorem~\ref{dPsurfaceintegral} still holds.  This is significant because classically, one is often interested in integral points on elliptic or hyperelliptic curves where the divisor ``at infinity'' is one or two points.

\end{document}